\documentclass[10pt, twoside]{article}
\usepackage{amsmath, amssymb, amsthm}
\long\def\symbolfootnote[#1]#2{\begingroup%
\def\thefootnote{\fnsymbol{footnote}}\footnote[#1]{#2}\endgroup}
\theoremstyle{plain}

\theoremstyle{definition}

\theoremstyle{remark}

\def\ntitle#1{\large \bf \noindent \begin{center} #1 \end{center}
\rm \normalsize }%
\def\nauthor#1{\noindent%
   \large \begin{center} \small \sc{#1} \end{center}
   \rm \vskip-13pt }
\let\title\ntitle
\let\author\nauthor

\let\address\naddress
\let\email\nemail

\begin{document}

\vspace{-.20in} \baselineskip .18in

\title{L-DIVERGENCE CONSISTENCY FOR A DISCRETE PRIOR}
\author{MARIAN GRENDAR} %% In regular caps as written.
\address{Department of Mathematics, FPV UMB  \\
 Tajovskeho 40, 974 01 Banska Bystrica, Slovakia\\
 Institute of Mathematics and CS, Banska Bystrica, Slovakia\\
 Institute of Measurement Sciences, Bratislava, Slovakia}
\email{marian.grendar@savba.sk} %%  your email address as written.
\centerline{\sc{summary}} \vspace{0.1in} \baselineskip.20in
\centerline{
\begin{minipage}{4.5in} \small{ \noindent
Posterior distribution over a countable set $\mathcal M$ of
continuous data-sampling
distributions %(DSD)
piles up at
$L$-projection of the true %DSD
distribution $r$ on $\mathcal M$, provided that the $L$-projection
is unique. If there are several $L$-projections of $r$ on $\mathcal
M$, then the posterior probability splits among them equally.
}\end{minipage}}

\vspace{.15in} \baselineskip.18in

\centerline{ \begin{minipage}{4.5in} \small{ \noindent {\it Keywords
and phrases:} Bayesian consistency, $L$-divergence, multiple
$L$-projections}\end{minipage}}

\vspace{.1in}

\centerline{ \begin{minipage}{4.5in} \small{ \noindent {\it AMS
Classification:} 60F10, 60F15}\end{minipage}}

\section{Introduction}

Walker \cite{W1} has recently considered consistency of posterior
distribution in Hellinger distance, for strictly positive prior over
a countable set of continuous data-sampling distributions. By means
of his martingale approach \cite{W2}, Walker developed a sufficient
condition for the Hellinger consistency of posterior density in the
above mentioned setting. Via a simple large-deviations approach we
show that in this setting posterior density is always consistent in
$L$-divergence. The consistency holds also under misspecification.
If there are multiple 'concentration points' ($L$-projections) the
posterior spreads among them equally.

\section{Bayesian nonparametric consistency}

Let there be countable set $\mathcal M = \{q_1, q_2, \dots \}$ of
probability density functions with respect to the Lebesgue measure;
sources, for short. On the set a Bayesian puts his strictly positive
prior probability mass function $\pi(\cdot)$. Let $r$ be the true
source of a random sample $X^n \triangleq X_1, X_2, \dots, X_n$.
Provided that $r \in \mathcal M$, as the sample size grows to
infinity, the posterior distribution $\pi(\cdot|X^n=x^n)$ over
$\mathcal M$ is expected to concentrate in a neighborhood of the
true source $r$. Whether and under what conditions this indeed
happens is a subject of Bayesian nonparametric consistency
investigations.  Surveys of the subject can be found at \cite{GGR},
\cite{WLP} among others.

Ghosal, Ghosh and Ramamoorthi \cite{GGR} define consistency of a
sequence of posteriors with respect to a metric or discrepancy
measure $d$ as follows: The sequence $\{\pi(\cdot|X^n), n \ge 1\}$
is said to be $d$-consistent at $r$, if there exists a $\Omega_0
\subset \mathbb R^\infty$ with $r(\Omega_0) = 1$ such that for
$\omega \in \Omega_0$, for every neighborhood $U$ of $r$,
$\pi(U|X^n) \rightarrow 1$ as $n$ goes to infinity. If a posterior
is $d$-consistent for any $r \in \mathcal M$ then it is said to be
$d$-consistent. There, two modes of convergence are usually
considered: convergence in probability and almost sure convergence.

{Obviously, in the definition the set of sources is not restricted
to be countable.} The present work is concerned with %investigates
the countable $\mathcal M$ case.

\section{Sanov's Theorem for Sources, L-consistency}

Let $\mathcal{M}^e \triangleq \{q: q \in \mathcal{M}, \pi(q) > 0\}$
be support of the prior pmf. In what follows, $r$ is not necessarily
from $\mathcal{M}^e$. Thus we are interested also in Bayesian
consistency under misspecification; i.e., when $\pi(r) = 0$. The
problem is the same as in the case of standard Bayesian consistency
(cf. Sect. 2): to find the source(s) upon which the posterior
concentrates.

For two densities $p, q$ with respect to the Lebesgue
measure\footnote{Any $\sigma$-finite measure, in general.}
$\lambda$, the $I$-divergence $I(p||q) \triangleq \int p \log(p/q)$.
The $L$-divergence $L(q||p)$ of $q$ with respect to $p$ is defined
as $L(q||p) \triangleq - \int p \log q$. The $L$-projection
$\hat q$ of $p$ on %set of sources
$\mathcal Q$ is $\hat q \triangleq \arg \inf_{q \in \mathcal Q}
L(q||p)$. There $\mathcal Q$ is a set of probability densities
defined on the same support. The value of $L$-divergence at an
$L$-projection of $p$ on $\mathcal Q$ is denoted by $L(\mathcal
Q||p)$.

The following  Sanov's Theorem for Sources ($L$ST) will be needed
for establishing the consistency in $L$-divergence. The Theorem
provides rate of the exponential decay of the posterior probability.

\smallskip

$L${\bf ST} \ \emph{Let $\mathcal N \subset \mathcal{M}^e$.} As $n
\rightarrow \infty$,
$$
\frac{1}{n} \log \pi(q \in \mathcal N | x^n) \rightarrow -
\{L(\mathcal N||r) - L(\mathcal{M}^e||r)\},
$$
\emph{with probability one}.

\smallskip

\emph{Proof} Let $l_n(q) \triangleq \exp({\sum_{l=1}^n \log
q(X_l)})$, $l_n(A) \triangleq \sum_{q \in A} l_n(q)$,  and
$\rho_n(q) \triangleq \pi(q) l_n(q)$, $\rho_n(A) \triangleq \sum_{q
\in A} \rho_n(q)$. In this notation $\pi(q \in \mathcal N|x^n) =
\frac{\rho_n(\mathcal N)}{\rho_n(\mathcal{M}^e)}$. The posterior
probability is bounded above and below as follows:
$$
\frac{\hat{\rho}_n({\mathcal N})}{\hat{l}_n(\mathcal{M}^e)}  \le
\pi(q \in \mathcal N|x^n) \le \frac{\hat{l}_n(\mathcal
N)}{\hat{\rho}_n({\mathcal{M}^e})},
$$
where $\hat{l}_n(A) \triangleq \sup_{q \in A} l_n(q)$,
$\hat{\rho}_n(A) \triangleq \sup_{q \in A} \rho_n(q)$.

$\frac{1}{n}(\log\hat{l}_n(\mathcal{N}) -
\log\hat{\rho}_n(\mathcal{M}^e))$ converges with probability one to
$L(\mathcal{M}^e||r) - L(\mathcal{N}||r)$. The same is the 'point'
of a.s. convergence of $\frac{1}{n}\log$ of the lower bound. $\qed$

\smallskip

Let for $\epsilon > 0$, $\mathcal{N}_\epsilon^C(\mathcal{M}^e)
\triangleq \{q: L(q||r) - L(\mathcal{M}^e||r) > \epsilon, q \in
\mathcal{M}^e\}$. Let $\mathcal{N}_\epsilon(\mathcal{M}^e)
\triangleq \mathcal{M}^e\backslash\mathcal{N}_\epsilon^C$.

\smallskip

{\bf Corollary} \emph{Let there be a finite  number of
$L$-projections of $r$ on $\mathcal{M}^e$. As $n \rightarrow
\infty$, $\pi(q \in \mathcal{N}_\epsilon^C(\mathcal{M}^e)| x^n)
\rightarrow 0$, with probability one.}

\smallskip

Standard Bayesian consistency follows as a special  $\pi(r) > 0$
case of the Corollary.

\smallskip

\section{Posterior Equi-concentration of Sources}

If there is more than one $L$-projection of $r$ on $\mathcal{M}^e$,
how is the posterior probability asymptotically spread among them?
This issue is 'in probability' answered by the next Theorem. Let
$\mathcal{N}^1_\epsilon \subset \mathcal{N}_\epsilon(\mathcal{M}^e)$
 contain (among other sources) just one $L$-projection of
$r$ on $\mathcal{M}^e$.

\smallskip

{\bf Theorem} \emph{Let there be $\mathrm k$ $L$-projections of $r$
on $\mathcal{M}^e$. Then for $n$ going to infinity,  $\pi(q \in
\mathcal{N}^1_\epsilon|x^n) \rightarrow \frac{1}{\mathrm k}$, in
probability.}

\smallskip

\emph{Proof} For any $\epsilon > 0$, there exists such $n_0$ that
for $n > n_0$, $r\{x^n: S(\hat{q}_\lambda) = S(\hat{q}_L)\} = 1$,
where $\hat{q}_\lambda \triangleq \arg \sup_{q \in \mathcal{M}^e}
\pi(q|x^n)$, $\hat{q}_L$ is $L$-projection of $r$ on
$\mathcal{M}^e$, and $S(\cdot)$ stands for 'set of all'.
Consequently, $\pi(\hat{q}_L|x^n) \ge \pi(q|x^n)$ for all $q \in
\mathcal{M}^e$. Posterior $\pi(q \in \mathcal{N}^1_\epsilon|x^n)$
can be expressed as $(1-A)/(\mathrm{k}(1-B))$, where $A \triangleq
\sum_{\sigma_1} \pi(q|x^n)/\pi(\hat{q}_L|x^n)$, $B \triangleq
\sum_{\sigma_2} \pi(q|x^n)/\mathrm{k}\pi(q|x^n)$; $\sigma_1
\triangleq \mathcal{N}_\epsilon^1\backslash\hat{q}_L$, $\sigma_2
\triangleq \mathcal{M}^e\backslash\bigcup_{j=1}^k\hat{q}_L^j$.
Markov's inequality implies that $\pi(q|x^n)/\pi(\hat{q}_L|x^n)$
converges to zero, in probability. Slutsky's Theorem then implies
that $A$, $B$ converges to zero, in probability. $\qed$

\section{EndNotes}

In order to place this note in context let us make a few comments.

1) An inverse of Sanov's Theorem has been established by Ganesh and
O'Con\-nell \cite{GO} for the case of sources with finite alphabet,
by means of formal large-deviations approach. Unaware of their work,
the present author developed in \cite{g} an inverse of Sanov's
Theorem for $n$-sources, for both discrete and continuous alphabet
and applied it to conditioning by rare sources problem and criterion
choice problem; cf. also \cite{GJ}.

2) At \cite{g} the concepts of $L$-divergence and $L$-projection
were introduced. See \cite{g} for a short discussion on why or why
not the 'new' divergence.

3) The present form of Sanov's Theorem for Sources ($L$ST) as well
as its proof are new.

4) Bayesian consistency under misspecification has already been
studied by Kleijn and van der Vaart \cite{KV} for general setting of
continuous prior on a set of continuous sources, using a  different
technique. The authors developed {\it sufficient} conditions for
somewhat related consistency (cf. Corollary 2.1 and Lemma 6.4 of
\cite{KV}) as well as rates of convergence. The equi-concentration
was not considered there.

\subsection{Acknowledgements}

Supported by VEGA grant 1/3016/06.

\bigskip

With a typo in statement and proof of $L$ST ($\mathcal{M}^e$ and
$\mathcal N$ were interchanged) this note appeared as: M. Grendar,
L-divergence consistency for a discrete prior, \emph{J. Stat. Res.},
40(1), 73-76, 2006.

\end{document}